# A note on uniqueness for linear evolution PDEs posed on the quarter-plane

Andreas Chatziafratis [1,3,*] and Spyridon Kamvissis [2,3]

**Abstract.** In this paper, we announce a rigorous approach to establishing uniqueness results, under certain conditions, for initial-boundary-value problems for a class of linear evolution partial differential equations (PDEs) formulated in a quarter-plane. We also effectively propose an algorithm for constructing non-uniqueness counter-examples which do not satisfy the said conditions. Our approach relies crucially on the rigorous analysis of regularity and asymptotic properties of integral representations derived formally via the celebrated Unified Transform Method for each such PDE. For uniqueness, this boundary behavior analysis allows for a careful implementation of an 'energy-estimate' argument on the semi-unbounded domain. Our ideas are elucidated via application of the present technique to two concrete examples, namely the heat equation and the linear KdV equation with Dirichlet data on the positive quadrant, under a particular set of conditions at the domain boundary and at infinity. Importantly, this is facilitated by delicate refinement of previous results concerning the boundary behavior analysis of these two celebrated models. In addition, we announce a uniqueness theorem for the linearized BBM equation, whose proof, in a similar spirit, will appear in a forthcoming paper. Finally, we briefly demonstrate how the general case of *oblique* Robin data can be recast as a Dirichlet problem. To the best of our knowledge, such well-posedness results appear for the first time in either the classical or the weak sense. Extensions to other classes of equations are underway and will appear elsewhere.

## 1. Introduction

About 25 years ago, Fokas introduced a technique [1], nowadays broadly known as Fokas method or unified transform method (UTM), for solving *integrable* nonlinear PDEs which was based on the synergy between the *Riemann-Hilbert* theory and the *Lax-pairs* formalism. It was soon realized [2-4] that this novel method offered a convenient and highly advantageous (in contrast to conventional transform methods and their well-known pathologies [4-8], even in the limited cases when the latter are 'applicable') framework for studying linear PDEs posed in semi-unbounded domains too. In the latter case, this novel method yields solution representations expressed as contour integrals in the spectral complex plane, and is thus regarded as the analogue of the *Fourier* transform which is traditionally involved in the solution of the respective problems on the whole space.

More specifically, as far as initial-boundary-value problems (IBVPs) for linear PDEs are concerned, the UTM methodology is a groundbreaking algorithmic procedure for formally deriving formulae (as *candidate* solutions) for evolution equations posed on a semi-infinite line or a finite interval, assuming existence and certain properties for the sought-after solution which allow the proposed steps of the recipe (as is always the case with transform methods), see e.g. [5,6]. In particular, the UTM consists of three main steps: (1) Formal construction, assuming validity of all calculations involved, of an integral representation, in the complex plane, for the unknown solution; this representation is not yet effective because it contains certain transforms of unknown boundary values. (2) Analysis of the so-called *global* relation, which is a simple equation connecting the given initial and boundary data with the unknown boundary values. (3) Utilization of certain invariance properties of the global relation dictated by the so-called *symmetry* relation combined with algebraic manipulations in order to eliminate the transforms of the unknown boundary values arising in the first step of the algorithm.

This modern method has been successful in a spectacularly broad class of linear-PDE problems (of which the greatest part is *not* amenable to treatment via classical methods) and has therefore been established as a highly efficient approach, both analytically and numerically, for applications in connection with areas as diverse as boundary value problems, spectral theory, control theory, inverse problems, and so on. For all these areas of application of the UTM, and extensions thereof by many mathematicians, see, for example, [9-41] and many more references cited therein. Incidentally, the UTM has as well led to the development of a new effective approach (pioneered by Fokas and Himonas who conceived and implemented it) to investigating well-posedness

---

[1] Department of Mathematics, National and Kapodistrian University of Athens

[2] Department of Pure and Applied Mathematics, University of Crete

[3] Institute of Applied and Computational Mathematics, FORTH, Crete

*Acknowledgements.*
A.C. thanks Professors: P.Biler, G. Biondini, J.L. Bona, J. Colliander, B. Deconinck, A.S. Fokas, M. Grillakis, T. Hatziafratis, A.A. Himonas, P. Kevrekidis, J. Lenells, R.L. Pego, B. Pelloni, D.A. Smith, I.G. Stratis and R.M. Temam, for providing encouragement and inspiration.

*corresponding author, e-mail: chatziafrati@math.uoa.gr

questions for nonlinear analogues, e.g. [42-47]; see also [48] for the recent extension of this approach by Himonas and Yan to *systems* of nonlinear PDEs.

However, surprisingly, the literature severely lacks complete and accurate results on existence and uniqueness of classical solutions to IBVPs for linear PDEs. Indicatively, in particular, as mentioned above, the derivation of a (candidate) solution formula for a given IBVP for a PDE, e.g. with *polynomial* dispersion relation, is vitally based on formal application of Green's theorem to the divergence form of the equation on a semi-infinite strip, which in turn a priori assumes, inter alia, sufficient regularity of the yet unknown solution as well as uniform –in time– decay for the solution and certain derivatives as $x \to +\infty$. Remarkably, the formulae constructed via the UTM were never (until the recent appearance [49] of a new line of rigorous investigations, see below) verified *a posteriori* (which would also guarantee *existence* of solutions). Apparently, this is a common pitfall in the literature dealing with 'exact solutions' via transform methods (classical or not). Moreover, in general, one cannot know whether the formula arrived at is a *unique* solution, since the derivation process is admittedly formal. Stated differently, even if the derivation were *rigorous* it would still yield *all candidate* explicit solutions; then, rigorous *back-substitution* (or justification) would indeed secure the validity of the proposed closed-form solution. Properly carrying out these fundamental and indispensable steps should have been the foremost task when invoking "the solution" to a specific model – otherwise one just risks merely doing calculus on a multi-variable function which ultimately may in general be completely un-related to the model under consideration.

In view of remedying the aforementioned careless reasoning throughout the literature, a new analytical approach was recently introduced in [49], applied in [50,51] and extended in, e.g., [52-60], for the rigorous refinement and generalization of the UTM as well as for analytical investigation of a miscellany of qualitative properties of such PDEs, e.g. constructive existence and spatiotemporal asymptotics. In the particular case of evolution PDEs with polynomial dispersion relation, the novel approach utilizes, as a starting point, the formulae afforded by suitable implementation of the UTM. One of the salient aspects of these works is the proper *a posteriori* justification of the validity of formally derived UTM representations, including reconstruction –in the sense of limits near the semi-infinite boundaries– of prescribed initial and boundary data – far from a trivial task, in the majority of cases. These rigorous investigations around the UTM have led to discovery of new instability, blow-up and break-down phenomena (e.g., [54-56] and forthcoming papers). Moreover, they have already had an impact on the works of several researchers in the area; for instance, compare the older [61-63] with [64-66], etc.

The purpose of this note is to start settling, at last, the (non)uniqueness issue. As additional motivation, we begin by presenting two *new* examples of *non*-uniqueness.

***Example 1*** The UTM leads to the following solution:

$$v(x,t) = \frac{i}{\pi}\int_{\gamma}[e^{i\lambda x-\lambda^2 t} - e^{i\lambda x}]\frac{d\lambda}{\lambda}, \text{ defined for } x>0,\ t>0,$$

to the IBVP: $v_t = v_{xx}$, $v(0,t) = \lim_{x\to 0^+} v(x,t) = 1$, $v(x,0) = \lim_{t\to 0^+} v(x,t) = 0$. (See [49,50].) The contour $\gamma$ is the oriented boundary of the domain $\{\lambda \in \mathbb{C} : \mathrm{Im}\lambda \ge 0, \mathrm{Re}(\lambda^2) \le 0\}$.

By Cauchy's theorem and Jordan's lemma, this solution can be written also in the following way:

$$v(x,t) = \frac{i}{\pi}\int_{-\infty}^{\infty}[e^{i\lambda x-\lambda^2 t} - e^{i\lambda x}]\frac{d\lambda}{\lambda} = \frac{i}{\pi}\int_{\gamma_0} e^{i\lambda x-\lambda^2 t}\frac{d\lambda}{\lambda},$$

where $\gamma_0 := (\gamma \cap \{|\lambda| \ge 1\}) + (\{\mathrm{Re}(\lambda^2) \ge 0\} \cap \kappa^+)$ and $\kappa^+$ is the semicircle $\kappa^+(s) = e^{i(\pi - s)}$, $0 \le s \le \pi$ (see Fig.1 and Fig.4).

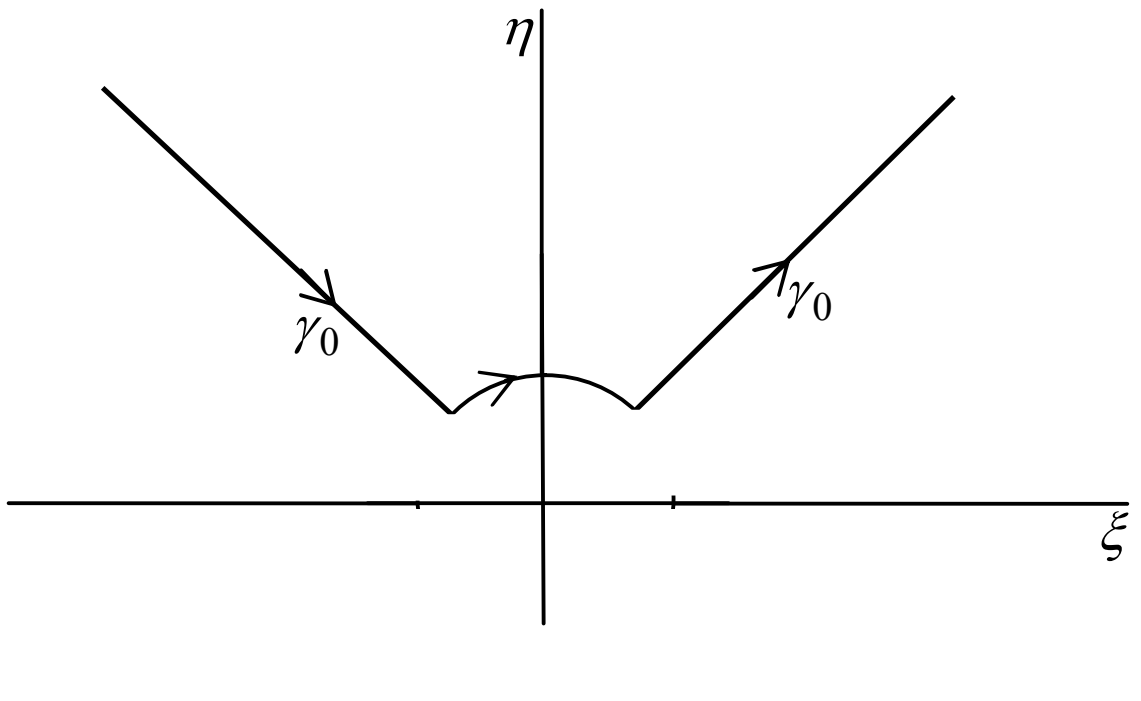

**Fig.1** *The contour $\gamma_0$ in the $\lambda-$plane ($\lambda = \xi + i\eta$)*

Its derivative with respect to $t$,

$$u(x,t)=\frac{\partial v(x,t)}{\partial t}=-\frac{i}{\pi}\int_{\gamma_0} e^{i\lambda x-\lambda^2 t}\lambda d\lambda=-\frac{i}{\pi}\int_{-\infty}^{\infty} e^{i\lambda x-\lambda^2 t}\lambda d\lambda=\frac{1}{2\sqrt{\pi}}\frac{x}{t\sqrt{t}}e^{-x^2/4t},$$

(which – upto a constant – is the $x$-derivative of Gauss's kernel), solves the IBVP:

$$u_t=u_{xx},\ u(0,t)=\lim_{x\to 0^+}u(x,t)=0,\ u(x,0)=\lim_{t\to 0^+}u(x,t)=0, \tag{1.1}$$

and $u(x,t)\not\equiv 0$.

Thus, the problem (1.1) does not have unique solution. As a matter of fact, each of the functions

$$u_n(x,t):=\frac{\partial^n v(x,t)}{\partial t^n}=-\frac{i}{\pi}\int_{-\infty}^{\infty}(-\lambda^2)^{n-1}e^{i\lambda x-\lambda^2 t}\lambda d\lambda=\frac{x}{2\sqrt{\pi}}\frac{\partial^{n-1}}{\partial t^{n-1}}\left(\frac{1}{t\sqrt{t}}e^{-x^2/4t}\right),\ n=1,2,3,\ldots,$$

solves (1.1).

***Example 2*** The UTM gives the following solution

$$v(x,t)=-\frac{3}{2\pi i}\int_{\operatorname{Im}\lambda=\varepsilon} e^{i\lambda x+i\lambda^3 t}\frac{d\lambda}{\lambda},\ \text{defined for } x>0,\ t>0,$$

to the IBVP: $v_t=-v_{xxx}$, $v(0,t)=\lim_{x\to 0^+}v(x,t)=1$, $v(x,0)=\lim_{t\to 0^+}v(x,t)=0$. (See [53].) The above integral does not depend on $\varepsilon$, where $\varepsilon$ is any fixed positive number.

Its derivative with respect to $t$,

$$u(x,t)=\frac{\partial v(x,t)}{\partial t}=-\frac{9}{2\pi}\int_{\operatorname{Im}\lambda=\varepsilon} e^{i\lambda x+i\lambda^3 t}\lambda^2 d\lambda=\frac{9}{2\pi}\frac{\partial^2}{\partial x^2}\left(\int_{\operatorname{Im}\lambda=\varepsilon} e^{i\lambda x+i\lambda^3 t}d\lambda\right), \tag{1.2}$$

solves the IBVP:

$$u_t=-u_{xxx},\ u(0,t)=\lim_{x\to 0^+}u(x,t)=0,\ u(x,0)=\lim_{t\to 0^+}u(x,t)=0. \tag{1.3}$$

Also we claim that $u(x,t)\not\equiv 0$. Indeed, if $u(x,t)\equiv 0$ then $\frac{\partial v(x,t)}{\partial t}\equiv 0$, and this would imply $v(x,t)\equiv v(x,0)\equiv 0$, which would contradict $v(0,t)=1$.

Also, each of the functions

$$u_n(x,t):=\frac{\partial^n v(x,t)}{\partial t^n}=-\frac{9}{2\pi}\int_{\operatorname{Im}\lambda=\varepsilon}(i\lambda^3)^{n-1}e^{i\lambda x+i\lambda^3 t}\lambda^2 d\lambda,\ n=1,2,3,\ldots,$$

solves (1.3).

The construction and justification of the above counterexamples to uniqueness rely vitally on the boundary behavior analysis provided in [50] and [53], respectively. In addition, these examples clearly indicate a novel 'recipe' for constructing such examples (not to be found in the existing literature) for a wide class of equations containing also higher-order evolution PDEs as well as PDEs involving lower-order terms too.

Now, we consider the following cases of IBVPs so as to illustrate our precise uniqueness considerations:

***Problem 1*** *Solve*

$$\begin{cases}\dfrac{\partial U}{\partial t}+\dfrac{\partial^3 U}{\partial x^3}=f,\ (x,t)\in Q:=\mathbb{R}^+\times\mathbb{R}^+\\ \lim\limits_{t\to 0^+}U(x,t)=u_0(x),\ x\in\mathbb{R}^+\\ \lim\limits_{x\to 0^+}U(x,t)=g_0(t),\ t\in\mathbb{R}^+,\end{cases} \tag{1.4}$$

*for* $U=U(x,t)$.

***Assumptions*** *Throughout this paper, we assume that*

$$u_0(x)\in\mathcal{S}([0,\infty)),\ g_0(t)\in C^\infty([0,\infty))\ \textit{and}\ f=f(x,t)\in C^\infty(\overline{Q}) \tag{1.5}$$

*such that* $f(\cdot,t)\in\mathcal{S}([0,\infty))$*, with respect to* $x$*, uniformy for* $t$ *in compact subsets of* $[0,+\infty)$.

***The UTM solution of Problem* 1 [6,9,53]** Setting

$$\hat{u}_0(\lambda) = \int_{y=0}^{\infty} e^{-i\lambda y} u_0(y)dy \text{ and } \hat{f}(\lambda,t) = \int_{y=0}^{\infty} e^{-i\lambda y} f(y,t)dy \text{, both defined for } \lambda \in \mathbb{C} \text{ with } \operatorname{Im}\lambda \le 0,$$

$$\tilde{g}_0(\omega(\lambda),t) = \int_{\tau=0}^{t} e^{\omega(\lambda)\tau} g_0(\tau)d\tau \text{, where } \omega(\lambda) = -i\lambda^3 \ (\lambda \in \mathbb{C}),$$

and

$$\tilde{\hat{f}}(\lambda,\omega(\lambda),t) = \int_{\tau=0}^{t} e^{\omega(\lambda)\tau} \hat{f}(\lambda,\tau)d\tau \ (\lambda \in \mathbb{C} \text{ with } \operatorname{Im}\lambda \le 0),$$

we define

$$\mathbb{I}_{\mathbb{R}}(x,t) = \int_{\lambda=-\infty}^{\infty} e^{i\lambda x - \omega(\lambda)t} \hat{u}_0(\lambda)d\lambda,$$

$$\mathbb{I}_{\Gamma}(x,t) = \int_{\lambda\in\Gamma} e^{i\lambda x - \omega(\lambda)t} [\alpha\hat{u}_0(\alpha\lambda) + \alpha^2\hat{u}_0(\alpha^2\lambda)]d\lambda,$$

$$\mathbb{I}_{g}(x,t) = \int_{\lambda\in\Gamma} e^{i\lambda x - \omega(\lambda)t} 3\lambda^2 \tilde{g}_0(\omega(\lambda),t)d\lambda,$$

$$\Phi_{\mathbb{R}}(x,t) = \int_{\lambda=-\infty}^{\infty} e^{i\lambda x - \omega(\lambda)t} \tilde{\hat{f}}(\lambda,\omega(\lambda),t)d\lambda,$$

$$\Phi_{\Gamma}(x,t) = \int_{\lambda\in\Gamma} e^{i\lambda x - \omega(\lambda)t} [\alpha\tilde{\hat{f}}(a\lambda,\omega(\lambda),t) + \alpha^2\tilde{\hat{f}}(\alpha^2\lambda,\omega(\lambda),t)]d\lambda,$$

where $\alpha = e^{2\pi i/3}$ and $\Gamma = \partial\Omega^-$ with $\Omega^- = \{\lambda \in \mathbb{C} : \operatorname{Im}\lambda \ge 0 \text{ and } \operatorname{Re}\omega(\lambda) \le 0\}$ (see Fig.2).

With this notation, the unified transform method gives the solution to the above problem (1.4) in the following form: For $x > 0$ and $t > 0$,

$$2\pi U(x,t) = \mathbb{I}_{\mathbb{R}}(x,t) + \mathbb{I}_{\Gamma}(x,t) - \mathbb{I}_{g}(x,t) + \Phi_{\mathbb{R}}(x,t) + \Phi_{\Gamma}(x,t). \tag{1.6}$$

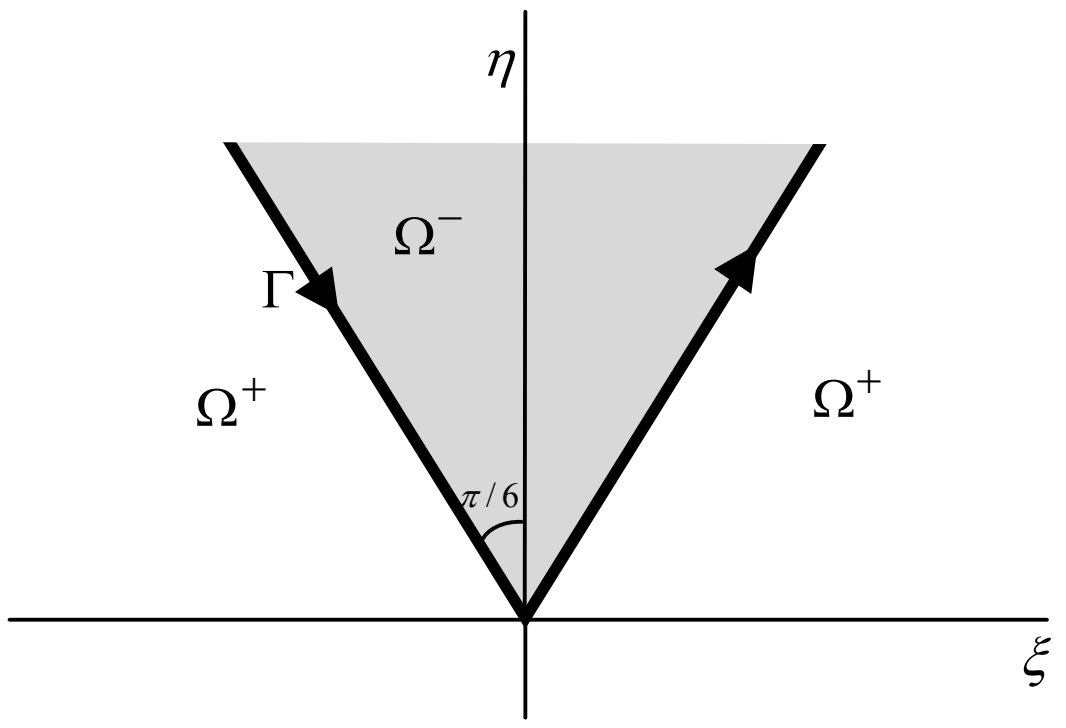

**Fig.2** *The contour* $\Gamma$ *is the boundary of* $\Omega^-$

Writing $\lambda = \xi + i\eta$, $\xi,\eta \in \mathbb{R}$, we have $\omega(\lambda) = -i(\xi + i\eta)^3 = (3\xi^2\eta - \eta^3) + i(-\xi^3 + 3\xi\eta^2)$.
Thus, $\Omega^- = \{\lambda \in \mathbb{C} : \eta \ge \sqrt{3}|\xi|\}$, $\Gamma = \{\lambda \in \mathbb{C} : \eta = \sqrt{3}|\xi|\}$. Also

$$\Omega^+ := \{\lambda \in \mathbb{C} : \operatorname{Im}\lambda \ge 0 \text{ and } \operatorname{Re}\omega(\lambda) \ge 0\} = \{\lambda \in \mathbb{C} : 0 \le \eta \le \sqrt{3}|\xi|\}.$$

In [53], we gave a rigorous proof that (1.6) is indeed a solution of (1.4) and we studied its boundary behavior, proving that the function $U(x,t)$, extended to $\overline{Q} - \{(0,0)\}$ by setting $U(x,0) = u_0(x)$ for $x > 0$ and $U(0,t) = g_0(t)$ for $t > 0$, is $C^\infty$ (in $\overline{Q} - \{(0,0)\}$).
We also proved that, for every nonnegative integers $k$, $m$, $\ell$, and for every $0 < t_0 < T_0$,

$$\lim_{x\to+\infty}\left( x^{\ell} \frac{\partial^{k+m}U(x,t)}{\partial x^k \partial t^m} \right) = 0, \tag{1.7}$$

uniformly for $t_0 \le t \le T_0$.
In this note, we will first prove that (1.7) is actually uniform for $0 < t \le T_0$. This refinement is essential in proving a uniqueness theorem for the solution (1.6) of equation (1.4).
Let us point out that (1.7) holds, uniformly for $0 < t \le T_0$, if and only if

$$\sup\left\{\left|x^{\ell}\frac{\partial^{k+m}U(x,t)}{\partial x^{k}\partial t^{m}}\right| : x>0,\, 0<t\le T_0\right\}<+\infty, \text{ for every nonnegative integers } k, m, \ell. \tag{1.8}$$

**Theorem 1** *The solution (1.6) of (1.4) satisfies (1.7), uniformly for* $0 < t \le T_0$ *(*$\forall T_0 > 0$*).*

**Theorem 2** *The function* $U(x,t)$*, defined by (1.6), is the unique solution of (1.4), in the following sense: If, in addition to (1.5), we assume that* $u_0(0) = g_0(0)$ *and* $g_0'(0) = -u_0'''(0) + f(0,0)$*, and*

$$V(x,t) \text{ is a } C^3 \text{ function in } \overline{Q} - \{(0,0\} \text{ and solves (1.4),} \tag{1.9}$$

$$\lim_{x\to\infty} V(x,t) = \lim_{x\to\infty} V_x(x,t) = 0 \text{ and } \sup_{x\ge 1}|V_{xx}(x,t)| < \infty \ (\forall t>0), \tag{1.10}$$

*and, for every* $T > 0$*, the functions*

$|V(x,t)|^2$, $|V_t(x,t)|^2$ *are, uniformly for* $0 < t \le T$*, integrable with respect to* $x \in [0,\infty)$*,*

*i.e., there exists a positive function* $\mathrm{B}_T(x)$ *such that* $\int_0^\infty \mathrm{B}_T(x)dx < +\infty$ *and, for* $0 < t \le T$*,*

$$|V(x,t)|^2 \le \mathrm{B}_T(x) \text{ and } |V_t(x,t)|^2 \le \mathrm{B}_T(x) \ (x>0), \tag{1.11}$$

*then* $V \equiv U$*.*

***Remarks***:

(1) For an example of a rather simplistic (and gravely inaccurate) formulation of a 'uniqueness theorem', the readers are referred to the Appendix of Mantzavinos *et al* [66]. This is a manifestation (thereby a reminder too) of the adage "the devil is in the details"!

(2) It should be highlighted that it is essential that one check that the candidate solution provided by the UTM indeed satisfies all conditions and hence that the set of solutions is not empty.

## 2. Proofs of Theorems 1 and 2

**2.1** ***Proof of Theorem 1 Step 1*** We claim that each of the functions $\mathbb{I}_\Gamma(x,t)$, $\mathbb{I}_g(x,t)$, $\Phi_\Gamma(x,t)$, satisfies (1.4), uniformly for $0 < t \le T_0$.
First, let us keep in mind that

$$\left|e^{i\lambda x}\right| = e^{-[x\sin(\pi/3)]|\lambda|} \text{ and } x^{\ell}\left|\lambda^N e^{i\lambda x}\right| \le x^{\ell}|\lambda|^N e^{-[x\sin(\pi/3)]|\lambda|}, \text{ for } \lambda\in\Gamma. \tag{2.1}$$

This implies that differentiation of the integrals $\mathbb{I}_\Gamma(x,t)$, $\mathbb{I}_g(x,t)$, $\Phi_\Gamma(x,t)$, with respect to $x$ or $t$, interchanges with the integration over $\Gamma$.
Now, to prove the claim for $\mathbb{I}_g(x,t)$, let us write the integral, which defines $\mathbb{I}_g(x,t)$, as follows:

$$\int_\Gamma \cdots d\lambda = \int_{\Gamma\cap\{|\lambda|\ge 1\}} \cdots d\lambda + \int_{\Gamma\cap\{|\lambda|\le 1\}} \cdots d\lambda .$$

In order to deal with the integral over $\Gamma \cap \{|\lambda| \ge 1\}$, it suffices to use the inequality

$$\left|e^{-\omega(\lambda)t}\tilde{g}_0(\omega(\lambda),t)\right| \le \int_0^t |g_0(\tau)|d\tau, \text{ for } \lambda\in\Gamma,$$

in combination with (2.1), to obtain that, for $x \ge 1$,

$$x^{\ell}\left|\int_{\Gamma\cap\{|\lambda|\ge 1\}}\lambda^N e^{i\lambda x-\omega(\lambda)t}3\lambda^2\tilde{g}_0(\omega(\lambda),t)d\lambda\right|\le 3\left(\int_0^t|g_0(\tau)|d\tau\right)x^{\ell}e^{-[x\sin(\pi/3)]/2}\int_{\Gamma\cap\{|\lambda|\ge 1\}}|\lambda|^{N+2}e^{-[\sin(\pi/3)]|\lambda|/2}d|\lambda|.$$

In order to deal with the integral over $\Gamma\cap\{|\lambda|\le 1\}$, it suffices to integrate by parts several times, depending on $\ell$. For example, writing

$$x\int_{\lambda\in\Gamma\cap\{|\lambda|\le 1\}}\lambda^N e^{i\lambda x-\omega(\lambda)t}\lambda^2\tilde{g}_0(\omega(\lambda),t)d\lambda=\frac{1}{i}\int_{\lambda\in\Gamma\cap\{|\lambda|\le 1\}}\frac{d(e^{i\lambda x})}{d\lambda}\lambda^N e^{-\omega(\lambda)t}\lambda^2\tilde{g}_0(\omega(\lambda),t)d\lambda$$

$$=\frac{1}{i}[\lambda^N e^{i\lambda x-\omega(\lambda)t}\lambda^2\tilde{g}_0(\omega(\lambda),t)\Big|_{\lambda=e^{i\pi/3}}^{\lambda=e^{i2\pi/3}}-\frac{1}{i}\int_{\lambda\in\Gamma\cap\{|\lambda|\le 1\}}e^{i\lambda x}\frac{d}{d\lambda}[\lambda^N e^{-\omega(\lambda)t}\lambda^2\tilde{g}_0(\omega(\lambda),t)]d\lambda,$$

we obtain that

$$\lim_{x\to\infty}\left[x\int_{\lambda\in\Gamma\cap\{|\lambda|\le 1\}}\lambda^N e^{i\lambda x-\omega(\lambda)t}\lambda^2\tilde{g}_0(\omega(\lambda),t)d\lambda\right]=0,\text{ uniformly for } 0<t\le T_0.$$

The proof of the claim for the functions $\mathbb{I}_\Gamma(x,t)$ and $\Phi_\Gamma(x,t)$ is similar.

***Step 2*** We claim that the function $\mathbb{I}_{\mathbb{R}}(x,t)$ satisfies (1.4), uniformly for $0<t\le T_0$.
To prove this, let us write

$$\mathbb{I}_{\mathbb{R}}(x,t)=\int_{-1}^{1}e^{i\lambda x-\omega(\lambda)t}\hat{u}_0(\lambda)d\lambda+\left(\int_{-\infty}^{-1}+\int_1^{\infty}\right)e^{i\lambda x-\omega(\lambda)t}[\hat{u}_0(\lambda)-\sigma_{\mathrm{M}}(\lambda)]d\lambda$$
$$+\int_{[-1+\sqrt{3}i,-1]\cup[1,1+\sqrt{3}i]}e^{i\lambda x-\omega(\lambda)t}\sigma_{\mathrm{M}}(\lambda)d\lambda+\int_{\Gamma\cap\{|\lambda|\ge 2\}}e^{i\lambda x-\omega(\lambda)t}\sigma_{\mathrm{M}}(\lambda)d\lambda,\quad(2.2)$$

where

$$\sigma_{\mathrm{M}}(\lambda):=\sum_{j=1}^{\mathrm{M}}\frac{d^{j-1}u_0}{dx^{j-1}}(0)\frac{1}{(i\lambda)^j}\quad(\lambda\in\mathbb{C},\ \lambda\ne 0).$$

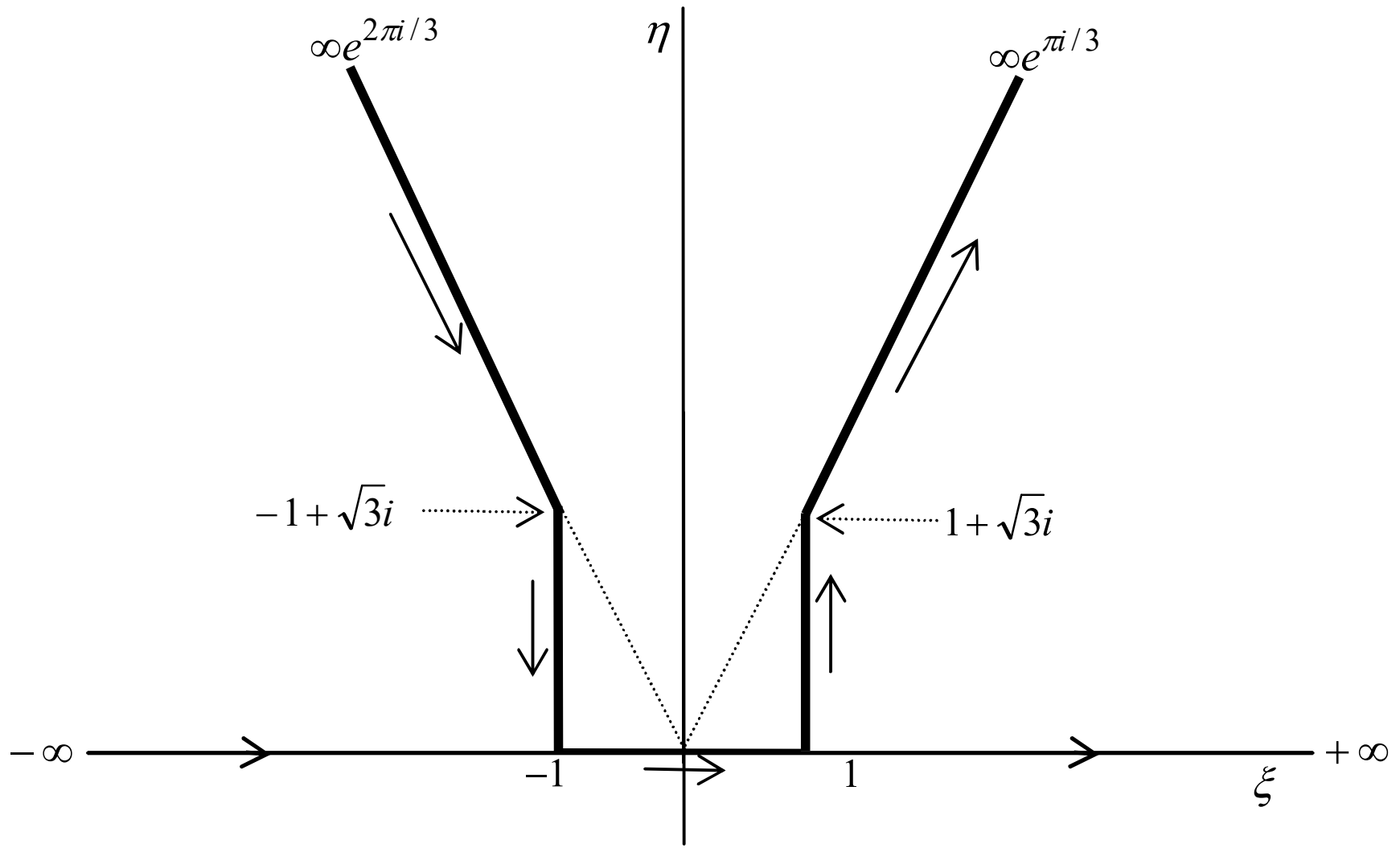

**Fig.3** *The contours of integration in (2.2)*

Keeping in mind that

$$\hat{u}_0(\lambda)-\sigma_{\mathrm{M}}(\lambda)=\frac{1}{(i\lambda)^{\mathrm{M}}}\int_{y=0}^{\infty}e^{-i\lambda y}\frac{d^{\mathrm{M}}u_0(y)}{dy^{\mathrm{M}}}dy=\mathrm{O}(1/\lambda^{\mathrm{M}+1}),\text{ as }\lambda\to\infty\text{ with }\operatorname{Im}\lambda\le 0,\quad(2.3)$$

and with $\mathrm{M}$ sufficiently large, we differentiate (2.2) and we obtain that, for $x>0$ and $t>0$,

$$\frac{\partial^{k+m}\mathbb{I}_{\mathbb{R}}(x,t)}{\partial x^k\partial t^m}=\int_{-1}^{1}(i\lambda)^k[-\omega(\lambda)]^m e^{i\lambda x-\omega(\lambda)t}\hat{u}_0(\lambda)d\lambda+\left(\int_{-\infty}^{-1}+\int_1^{\infty}\right)(i\lambda)^k[-\omega(\lambda)]^m e^{i\lambda x-\omega(\lambda)t}[\hat{u}_0(\lambda)-\sigma_{\mathrm{M}}(\lambda)]d\lambda$$

$$+\int_{[-1+\sqrt{3}i,-1]\cup[1,1+\sqrt{3}i]}(i\lambda)^k[-\omega(\lambda)]^m e^{i\lambda x-\omega(\lambda)t}\sigma_{\mathrm{M}}(\lambda)d\lambda+\int_{\Gamma\cap\{|\lambda|\geq 2\}}(i\lambda)^k[-\omega(\lambda)]^m e^{i\lambda x-\omega(\lambda)t}\sigma_{\mathrm{M}}(\lambda)d\lambda. \quad (2.4)$$

Multiplying (2.4) by $ix$ and integrating by parts we see that

$$(ix)\frac{\partial^{k+m}\mathbb{I}_{\mathbb{R}}(x,t)}{\partial x^k\partial t^m}=\int_{-1}^{1}e^{i\lambda x}\frac{d}{d\lambda}[(i\lambda)^k[-\omega(\lambda)]^m e^{-\omega(\lambda)t}\hat{u}_0(\lambda)]d\lambda$$

$$+\left(\int_{-\infty}^{-1}+\int_{1}^{\infty}\right)e^{i\lambda x}\frac{d}{d\lambda}\{(i\lambda)^k[-\omega(\lambda)]^m e^{-\omega(\lambda)t}[\hat{u}_0(\lambda)-\sigma_{\mathrm{M}}(\lambda)]\}d\lambda$$

$$+\int_{[-1+\sqrt{3}i,-1]\cup[1,1+\sqrt{3}i]}e^{i\lambda x}\frac{d}{d\lambda}\{(i\lambda)^k[-\omega(\lambda)]^m e^{-\omega(\lambda)t}\sigma_{\mathrm{M}}(\lambda)\}d\lambda$$

$$+\int_{\Gamma\cap\{|\lambda|\geq 2\}}e^{i\lambda x}\frac{d}{d\lambda}\{(i\lambda)^k[-\omega(\lambda)]^m e^{-\omega(\lambda)t}\sigma_{\mathrm{M}}(\lambda)\}d\lambda. \quad (2.5)$$

Let us point out that, in the above integration by parts processes, the boundary terms, coming from the intermediate evaluations, cancel each other.

Now, we can show, working as in *Step 1*, that the last integral in RHS of (2.5), which is taken on $\Gamma\cap\{|\lambda|\geq 2\}$, tends to zero, as $x\to+\infty$.

Also, setting $\lambda=1+si$, $0\leq s\leq\sqrt{3}$, the exponential $e^{i\lambda x}$ becomes

$$e^{i\lambda x}\Big|_{\lambda=1+si}=e^{ix}e^{-sx}\to 0,\ \text{as } x\to\infty,\ \text{for } s>0.$$

It follows, from Lebesgue's dominated convergence theorem, that the integral in (2.5), taken on the interval $[1,1+\sqrt{3}i]$, tends to zero, as $x\to\infty$. Similarly, we treat the integral taken on the interval $[-1+\sqrt{3}i,-1]$.

On the other hand, the other two integrals in (2.5), which are taken on subsets of the real line, also tend to zero, as $x\to+\infty$, by the Riemann-Lebesgue lemma.

Therefore,

$$\lim_{x\to\infty}\left[(ix)\frac{\partial^{k+m}\mathbb{I}_{\mathbb{R}}(x,t)}{\partial x^k\partial t^m}\right]=0,\ \text{uniformly for } 0<t\leq T_0.$$

Further integrations by parts show that

$$\lim_{x\to\infty}\left[(ix)^{\ell}\frac{\partial^{k+m}\mathbb{I}_{\mathbb{R}}(x,t)}{\partial x^k\partial t^m}\right]=0,\ \text{uniformly for } 0<t\leq T_0,\ \forall\ell,$$

and this completes the proof of the claim.

***Step 3*** We claim that the function $\Phi_{\mathbb{R}}(x,t)$ satisfies (1.4), uniformly for $0<t\leq T_0$.

As in the previous step, we have

$$\hat{f}(\lambda,t)=h_{\mathrm{M}}(\lambda,t)+\frac{1}{(i\lambda)^{\mathrm{M}}}\int_{y=0}^{\infty}e^{-i\lambda y}\frac{\partial^{\mathrm{M}}f(y,t)}{\partial y^{\mathrm{M}}}dy,\ \text{where } h_{\mathrm{M}}(\lambda,t):=\sum_{j=1}^{\mathrm{M}}\frac{1}{(i\lambda)^j}\frac{\partial^{j-1}f(y,t)}{\partial y^{j-1}}\bigg|_{y=0},$$

and

$$\hat{f}(\lambda,t)-h_{\mathrm{M}}(\lambda,t)=\mathrm{O}(1/\lambda^{\mathrm{M}+1}),\ \text{for } \lambda\to\infty,\ \text{with } \operatorname{Im}\lambda\leq 0,\ \text{uniformly for } 0<t\leq T_0. \quad (2.6)$$

Thus, we may write the integral $\Phi_{\mathbb{R}}(x,t)$ as follows: For any $\mathrm{M}\in\mathbb{N}$,

$$\Phi_{\mathbb{R}}(x,t)=\int_{-1}^{1}e^{i\lambda x-\omega(\lambda)t}\tilde{\hat{f}}(\lambda,\omega(\lambda),t)d\lambda+\left(\int_{-\infty}^{-1}+\int_{1}^{\infty}\right)e^{i\lambda x-\omega(\lambda)t}[\tilde{\hat{f}}(\lambda,\omega(\lambda),t)-\tilde{h}_{\mathrm{M}}(\lambda,\omega(\lambda),t)]d\lambda$$

$$+\int_{[-1+\sqrt{3}i,-1]\cup[1,1+\sqrt{3}i]}e^{i\lambda x-\omega(\lambda)t}\tilde{h}_{\mathrm{M}}(\lambda,\omega(\lambda),t)d\lambda+\int_{\Gamma\cap\{|\lambda|\geq 2\}}e^{i\lambda x-\omega(\lambda)t}\tilde{h}_{\mathrm{M}}(\lambda,\omega(\lambda),t)d\lambda. \quad (2.7)$$

In view of (2.6),

$$\tilde{\hat{f}}(\lambda,\omega(\lambda),t)-\tilde{h}_{\mathrm{M}}(\lambda,\omega(\lambda),t)=\mathrm{O}(1/\lambda^{\mathrm{M}+1}),\ \text{for } \mathbb{R}\ni\lambda\to\pm\infty,\ \text{uniformly for } 0<t\leq T_0. \quad (2.8)$$

Now, working as in *Step 2*, we can easily carry out the proof of the claim, using (2.7) and (2.8).

**2.2 *Proof of Theorem 2*** We may assume that the data $u_0$, $g_0$ and $f$, and the functions $U$ and $V$ are real-valued. Let us set

$$W(x,t) := V(x,t) - U(x,t).$$

Then

$$W(x,0) = 0 \text{ for } x > 0,\ W(0,t) = 0 \text{ for } t > 0,\ W_t = -W_{xxx} \text{ for } (x,t) \in \overline{Q} - \{(0,0\},$$

and (1.6), (1.7) and (1.8) hold with $V$ replaced by $U$ (by taking $M_T$ larger – if necessary). Indeed, this last assertion for $U$ follows from Theorem 1 and the Theorems in [53], using also the assumption $u_0(0) = g_0(0)$. Thus, (1.6), (1.7) and (1.8) hold also with $V$ replaced by $W$ (by taking $M_T$ larger – if necessary).

Then the equation $W_t = -W_{xxx}$ implies that

$$\int_{x=0}^{\infty} W(x,t) W_t(x,t)\,dx = -\int_{x=0}^{\infty} W(x,t) W_{xxx}(x,t)\,dx, \text{ for every } t > 0. \tag{2.9}$$

We note that, since (1.8) holds with $V$ replaced by $W$, the integrals in (2.9) are absolutely convergent. This also implies that, by Lebesgue's dominated convergence theorem,

$$\int_{x=0}^{\infty} W(x,t) W_t(x,t)\,dx = \frac{1}{2}\frac{d}{dt}\left[\int_{x=0}^{\infty} [W(x,t)]^2\,dx\right], \text{ for every } t > 0. \tag{2.10}$$

Integrating by parts, we obtain that, for $t > 0$ and $A > 0$,

$$\begin{aligned} -\int_{x=0}^{A} W(x,t) W_{xxx}(x,t)\,dx &= -[W(x,t) W_{xx}(x,t)]\Big|_{x=0}^{x=A} + \int_{x=0}^{A} W_x(x,t) W_{xx}(x,t)\,dx \\ &= -[W(x,t) W_{xx}(x,t)]\Big|_{x=0}^{x=A} + \frac{1}{2}\{[W_x(x,t)]^2\}\Big|_{x=0}^{x=A}. \end{aligned} \tag{2.11}$$

Letting $A \to \infty$ and taking into consideration that (1.7) holds with $V$ replaced by $W$, we see that (2.11) implies that

$$-\int_{x=0}^{\infty} W(x,t) W_{xxx}(x,t)\,dx = -\frac{1}{2}[W_x(0,t)]^2, \text{ for every } t > 0. \tag{2.12}$$

It follows from (2.9), (2.10) and (2.12), that

$$\frac{d}{dt}\left[\int_{x=0}^{\infty} [W(x,t)]^2\,dx\right] \le 0, \text{ for every } t > 0,$$

and, therefore,

$$\int_{x=0}^{\infty} [W(x,s)]^2\,dx \le \int_{x=0}^{\infty} [W(x,t)]^2\,dx, \text{ for every } s > t > 0. \tag{2.13}$$

Since

$$|V(x,t)|^2 \le M_T(x),\ 0 < t \le T, \text{ and } \int_0^{\infty} M_T(x)\,dx < +\infty,$$

Lebesgue's dominated convergence theorem implies that

$$\lim_{t \to 0^+} \int_{x=0}^{\infty} [W(x,t)]^2\,dx = \int_{x=0}^{\infty} [W(x,0)]^2\,dx = 0. \tag{2.14}$$

(The last equation in (2.14) follows from the fact that $W(x,0) = 0$, for $x > 0$.)

Thus, letting $t \to 0^+$ in (2.13), we obtain, in view of (2.14), that

$$\int_{x=0}^{\infty} [W(x,s)]^2\,dx \le 0\ (\forall s > 0) \Rightarrow W(x,s) \equiv 0.$$

This completes the proof.

## 3. A uniqueness theorem for the heat equation on the half-line

In this section we will study questions, similar to the ones in Sections 1 and 2, for the following problem.

***Problem 2*** *Assuming (1.5), solve*

$$\begin{cases} \dfrac{\partial U}{\partial t} - \dfrac{\partial^2 U}{\partial x^2} = f, \ (x,t) \in Q := \mathbb{R}^+ \times \mathbb{R}^+ \\ \lim\limits_{t \to 0^+} U(x,t) = u_0(x), \ x \in \mathbb{R}^+ \\ \lim\limits_{x \to 0^+} U(x,t) = g_0(t), \ t \in \mathbb{R}^+, \end{cases} \tag{3.1}$$

*for* $U = U(x,t)$.

**3.1 *The UTM solution for Problem* 2 [6,9,50]** For $x > 0$ and $t > 0$,

$$2\pi U(x,t) = \int_{-\infty}^{\infty} e^{i\lambda x - \lambda^2 t} \hat{u}_0(\lambda) d\lambda - \int_{\gamma} e^{i\lambda x - \lambda^2 t} \hat{u}_0(-\lambda) d\lambda - \int_{\gamma} e^{i\lambda x - \lambda^2 t} \lambda \tilde{g}_0(\lambda,t) d\lambda$$
$$+ \int_{-\infty}^{\infty} e^{i\lambda x - \lambda^2 t} \tilde{\hat{f}}(\lambda,t) d\lambda - \int_{\gamma} e^{i\lambda x - \lambda^2 t} \tilde{\hat{f}}(-\lambda,t) d\lambda, \tag{3.2}$$

where

$$\tilde{g}_0(\lambda,t) = \int_{\tau=0}^{t} e^{\lambda^2 \tau} g_0(\tau) d\tau \ (\lambda \in \mathbb{C}) \text{ and } \tilde{\hat{f}}(\lambda,t) = \int_{\tau=0}^{t} e^{\lambda^2 \tau} \hat{f}(\lambda,\tau) d\tau \ (\lambda \in \mathbb{C} \text{ with } \operatorname{Im}\lambda \le 0),$$

the contour $\gamma$ is the oriented boundary of the domain $\{\lambda \in \mathbb{C} : \operatorname{Im}\lambda \ge 0 \text{ } and \text{ } \operatorname{Re}(\lambda^2) \le 0\}$ (see Fig.4) and $\hat{u}_0(\lambda)$ and $\hat{f}(\lambda,\tau)$ are as in Section 1.

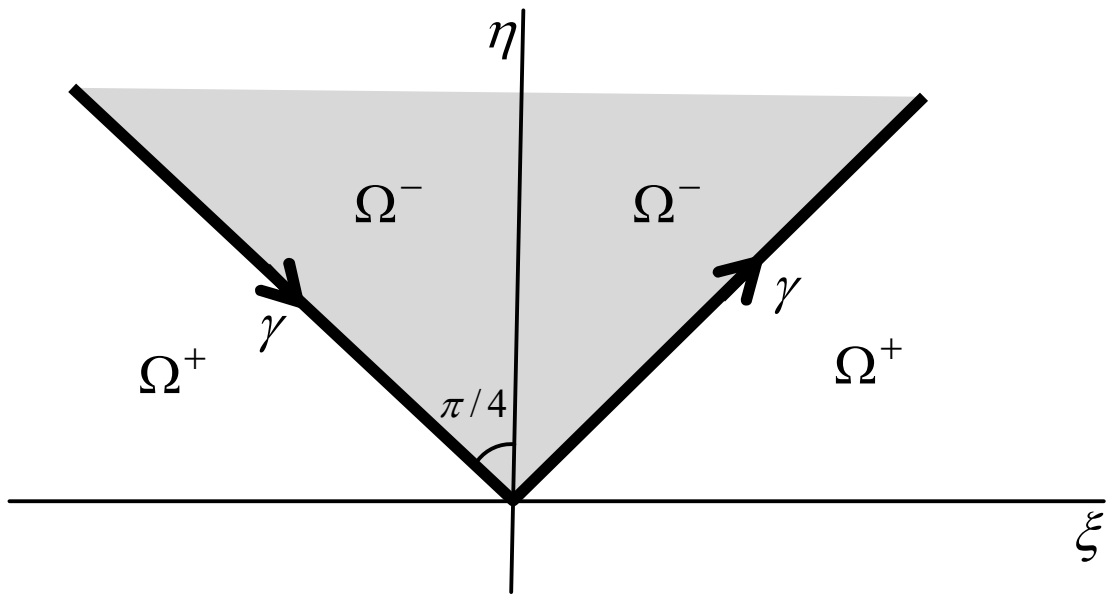

**Fig.4** *The contour* $\gamma = \partial(\{\lambda \in \mathbb{C} : \operatorname{Im}\lambda \ge 0 \text{ } and \text{ } \operatorname{Re}(\lambda^2) = 0\})$

**Theorem 3** *The solution* $U(x,t)$ *of (3.1), defined by (3.2), satisfies*

$$\lim_{x \to +\infty} \left( x^{\ell} \frac{\partial^{k+m} U(x,t)}{\partial x^k \partial t^m} \right) = 0, \tag{3.3}$$

*for every nonnegative integers* $k$, $m$, $\ell$, *uniformly for* $0 < t \le T_0$ *(*$\forall T_0 > 0$*).*

**Theorem 4** *The fuction* $U(x,t)$, *defined by (3.2), is the unique solution of (3.1), in the following sense: If, in addition to (1.2), we assume that* $u_0(0) = g_0(0)$ *and* $u_0''(0) + f(0,0) = g_0'(0)$, *and*

$$V(x,t) \text{ is a } C^2 \text{ function in } \overline{Q} - \{(0,0\} \text{ and solves (3.1)}, \tag{3.4}$$

$$\lim_{x \to \infty} V(x,t) = 0 \text{ and } \sup_{x \ge 1} |V_x(x,t)| < \infty \ (\forall t > 0), \tag{3.5}$$

*and, for every* $T > 0$,

$$|V(x,t)|^2, \ |V_t(x,t)|^2 \text{ are, uniformly for } 0 < t \le T, \text{ integrable with respect to } x \in [0,\infty),$$

*i.e., there exists a positive function* $\mathrm{B}_T(x)$ *such that* $\int_0^{\infty} \mathrm{B}_T(x) dx < +\infty$ *and, for* $0 < t \le T$,

$$|V(x,t)|^2 \le \mathrm{B}_T(x) \text{ and } |V_t(x,t)|^2 \le \mathrm{B}_T(x) \ (x > 0), \tag{3.6}$$

*then* $V \equiv U$.

**3.2 *Proof of Theorem 3*** Firstly, as in *Step 1* of the proof of Theorem 1, we can show that the functions defined by the integrals in (3.2), which are taken on $\gamma$, satisfy (3.3).

Now the first integral in the RHS of (3.2) can be written as follows:

$$\int_{-\infty}^{\infty} e^{i\lambda x-\lambda^2 t}\hat{u}_0(\lambda)d\lambda = \int_{-1}^{1} e^{i\lambda x-\lambda^2 t}\hat{u}_0(\lambda)d\lambda + \left(\int_{-\infty}^{-1} + \int_{1}^{\infty}\right) e^{i\lambda x-\lambda^2 t}[\hat{u}_0(\lambda)-\sigma_{\mathrm{M}}(\lambda)]d\lambda$$

$$+ \int_{[-1+i,-1]\cup[1,1+i]} e^{i\lambda x-\omega(\lambda)t}\sigma_{\mathrm{M}}(\lambda)d\lambda + \int_{\gamma\cap\{|\lambda|\geq\sqrt{2}\}} e^{i\lambda x-\omega(\lambda)t}\sigma_{\mathrm{M}}(\lambda)d\lambda, \quad (3.7)$$

where $\sigma_{\mathrm{M}}(\lambda)$ is as in Section 2. Choosing $\mathrm{M}$ sufficiently large, we see that the function defined by the integral in LHS of (3.7) satisfies (3.3).
Also, similarly to (2.7), we can write

$$\int_{-\infty}^{\infty} e^{i\lambda x-\lambda^2 t}\tilde{f}(\lambda,t)d\lambda = \int_{-1}^{1} e^{i\lambda x-\lambda^2 t}\tilde{f}(\lambda,t)d\lambda + \left(\int_{-\infty}^{-1} + \int_{1}^{\infty}\right) e^{i\lambda x-\lambda^2 t}[\tilde{f}(\lambda,t)-\tilde{h}_{\mathrm{M}}(\lambda,t)]d\lambda$$

$$+ \int_{[-1+i,-1]\cup[1,1+i]} e^{i\lambda x-\lambda^2 t}\tilde{h}_{\mathrm{M}}(\lambda,t)d\lambda + \int_{\gamma\cap\{|\lambda|\geq\sqrt{2}\}} e^{i\lambda x-\lambda^2 t}\tilde{h}_{\mathrm{M}}(\lambda,t)d\lambda, \quad (3.8)$$

where

$$\tilde{h}_{\mathrm{M}}(\lambda,t) = \int_{\tau=0}^{t} e^{\lambda^2\tau}h_{\mathrm{M}}(\lambda,\tau)d\tau \ (\lambda\in\mathbb{C}) \text{ and } h_{\mathrm{M}}(\lambda,\tau) \text{ is as in Section 2.}$$

The above conclusions easily imply that, indeed, $U(x,t)$ satisfies (3.3).

**3.3 *Proof of Theorem 4*** We may assume that the data $u_0$, $g_0$ and $f$, and the functions $U$ and $V$ are real-valued. Let us set

$$W(x,t) := V(x,t) - U(x,t).$$

Then

$$W(x,0)=0 \text{ for } x>0,\ W(0,t)=0 \text{ for } t>0,\ W_t = W_{xx} \text{ for } (x,t)\in\overline{Q}-\{(0,0\},$$

and (3.4), (3.5) and (3.6) hold with $V$ replaced by $U$ (by taking $M_T$ larger – if necessary). Indeed, this last assertion for $U$ follows from Theorem 3 and the Theorems in [11], using also the assumption $u_0(0)=g_0(0)$. Thus, (3.4), (3.5) and (3.6) hold also with $V$ replaced by $W$ (by taking $M_T$ larger – if necessary).
Then the equation $W_t = W_{xx}$ implies that

$$\int_{x=0}^{\infty} W(x,t)W_t(x,t)dx = \int_{x=0}^{\infty} W(x,t)W_{xx}(x,t)dx, \text{ for every } t>0. \quad (3.9)$$

We note that, since (3.6) holds with $V$ replaced by $W$, the integrals in (3.9) are absolutely convergent. This also implies that, by Lebesgue's dominated convergence theorem,

$$\int_{x=0}^{\infty} W(x,t)W_t(x,t)dx = \frac{1}{2}\frac{d}{dt}\left[\int_{x=0}^{\infty}[W(x,t)]^2dx\right], \text{ for every } t>0. \quad (3.10)$$

Integrating by parts, we obtain that, for $t>0$ and $A>0$,

$$\int_{x=0}^{A} W(x,t)W_{xx}(x,t)dx = [W(x,t)W_x(x,t)]\Big|_{x=0}^{x=A} - \int_{x=0}^{A}[W_x(x,t)]^2dx. \quad (3.11)$$

Letting $A\to\infty$ and taking into consideration that (3.5) holds with $V$ replaced by $W$, we see that (3.11) implies that

$$\int_{x=0}^{\infty} W(x,t)W_{xx}(x,t)dx = -\int_{x=0}^{A}[W_x(x,t)]^2dx, \text{ for every } t>0. \quad (3.12)$$

It follows from (3.9), (3.10) and (3.12), that

$$\frac{d}{dt}\left[\int_{x=0}^{\infty}[W(x,t)]^2dx\right] \leq 0, \text{ for every } t>0,$$

and, therefore,

$$\int_{x=0}^{\infty}[W(x,s)]^2dx \leq \int_{x=0}^{\infty}[W(x,t)]^2dx, \text{ for every } s>t>0. \quad (3.13)$$

Since

$$|V(x,t)|^2 \le M_T(x)\,,\ 0<t\le T\ \text{, and } \int_0^\infty M_T(x)dx < +\infty\,,$$

Lebesgue's dominated convergence theorem implies that

$$\lim_{t\to 0^+} \int_{x=0}^{\infty} [W(x,t)]^2 dx = \int_{x=0}^{\infty} [W(x,0)]^2 dx = 0\,. \tag{3.14}$$

(The last equation in (3.14) follows from the fact that $W(x,0)=0$, for $x>0$.)

Thus, letting $t\to 0^+$ in (3.13), we obtain, in view of (3.14), that

$$\int_{x=0}^{\infty} [W(x,s)]^2 dx \le 0\ \ (\forall s>0)\ \Rightarrow\ W(x,s)\equiv 0\,.$$

This completes the proof.

***Remarks:***

1) It is obvious that the above constructive technique for proving uniqueness is straightforwardly applicable to a very large class of evolution PDEs which may involve also higher- and/or lower- order terms with appropriate signs in their coefficients allowing for the implementation of the energy-argument step.

2) It may be checked that proving uniqueness of solution for IBVPs (for a given evolution PDE) with other types of mixed data reduces to the Dirichlet case, by use of an appropriate transformation. This point is elaborated in the next section.

## 4. Uniqueness theorems for the BBM equation and for Robin-type problems on the half-line

Let us first consider the following Dirichlet problem for the linearized Benjamin-Bona-Mahony (BBM) equation. An appropriate extension of the UTM for the solution (with rigorous verification followed by full boundary behavior analysis, well-posedness, spatiotemporal asymptotics, and so forth) of Problem 3, below, will appear in [59].

***Problem* 3** *Given $u_0(x)$, $g_0(t)$ and $f(x,t)$, solve*

$$\begin{cases} u_t - \alpha u_{xxt} + \beta u_x = f,\ (x,t)\in Q \\ u(x,0) = u_0(x),\ x\in\mathbb{R}^+ \\ u(0,t) = g_0(t),\ t\in\mathbb{R}^+, \end{cases} \tag{4.1}$$

*for $u=u(x,t)$ ($\alpha,\beta>0$).*

**Theorem 5 [59]** *Assuming (1.5) and $u_0(0)=g_0(0)$, the UTM solution $u$ of (4.1) is unique in the following sense: If $v(x,t)$ is $C^3$ in $\overline{Q}-\{(0,0)\}$, satisfies (4.1),*

$$\lim_{x\to\infty} v(x,t) = 0\,,\ \sup_{x\ge 1}|v_{xt}(x,t)| < \infty\ (\forall t>0), \tag{4.2}$$

*and, for every $T>0$, the functions $|v(x,t)|^2$, $|v_x(x,t)|^2$, $|v_t(x,t)|^2$, $|v_{xt}(x,t)|^2$, are, uniformly for $0<t\le T$, integrable with respect to $x\in[0,\infty)$, i.e., there exists a positive function $\mathcal{M}_T(x)$ such that $\int_0^\infty \mathcal{M}_T(x)dx < +\infty$ and, for $0<t\le T$ and $x>0$,*

$$|v(x,t)|^2 \le \mathcal{M}_T(x)\,,\ |v_x(x,t)|^2 \le \mathcal{M}_T(x)\,,\ |v_t(x,t)|^2 \le \mathcal{M}_T(x)\,,\ |v_{xt}(x,t)|^2 \le \mathcal{M}_T(x)\,, \tag{4.3}$$

*then $v\equiv u$.*

Next, we consider the following IBVP (also dealt with in [59]) with Robin data:

***Problem* 4** *Given* $u_0(x)$, $g_0(t)$ *and* $f(x,t)$, *solve*

$$\begin{cases} u_t - \alpha u_{xxt} + \beta u_x = f, \ (x,t) \in Q := \mathbb{R}^+ \times \mathbb{R}^+ \\ u(x,0) = u_0(x), \ x \in \mathbb{R}^+ \\ Au(0,t) + Bu_x(0,t) = g_0(t), \ t \in \mathbb{R}^+, \end{cases} \tag{4.4}$$

*for* $u = u(x,t)$.

***Reduction of uniqueness of Problem 4 to the uniqueness of Problem 3*** Let $u$ be the UTM solution of (4.4) [59] and $v$ is a solution of (4.4), satisfying appropriate conditions [59]. Setting $w = (Av + Bv_x) - (Au + Bu_x)$, we have that

$$w_t - w_{xxt} + w_x = 0, \ w(x,0) = 0 \text{ and } w(0,t) = 0.$$

By Theorem 5, we obtain that $w \equiv 0$, i.e., $A\xi + B\xi_x \equiv 0$, where we have set $\xi(x,t) = v(x,t) - u(x,t)$. Since the case $B = 0$ is covered by Theorem 5, we may assume that $B \neq 0$, in which case we see that the function $\zeta(x,t) =: \xi(x,t)\exp[Ax/B]$ satisfies the equation $\zeta_x \equiv 0$. Hence, $\zeta(x,t)$ is a function of $t$, let us say $\zeta(x,t) = \phi(t)$. Then $\xi(x,t) = \phi(t)\exp[-Ax/B]$ and $\phi(0) = 0$. Also

$$\xi_t - \xi_{xxt} + \xi_x = \phi'(t)e^{-Ax/B} - \frac{A^2}{B^2}\phi'(t)e^{-Ax/B} - \frac{A}{B}\phi(t)e^{-Ax/B} = 0 \ \Rightarrow \ \left(1 - \frac{A^2}{B^2}\right)\phi'(t) - \frac{A}{B}\phi(t) = 0.$$

Therefore $\phi \equiv 0$, i.e., $v \equiv u$.

We now conclude this paper by illustrating how an IBVP with *oblique* Robin conditions is also reducible to the corresponding Dirichlet-type IBVP for the same PDE.

***Problem* 5** *Given* $u_0(x)$, $g_0(t)$ *and* $f(x,t)$, *solve*

$$\begin{cases} u_t = u_{xx} + f, \ (x,t) \in Q := \mathbb{R}^+ \times \mathbb{R}^+ \\ u(x,0) = u_0(x), \ x \in \mathbb{R}^+ \\ Au(0,t) + Bu_x(0,t) + Cu_t(0,t) = g_0(t), \ t \in \mathbb{R}^+, \end{cases} \tag{4.5}$$

*for* $u = u(x,t)$ *(*$A, B, C > 0$*).*

(A formula for Problem 5 was proposed in [61,62]; a meticulous analysis about this problem as well as for the linearized Boussinesq equation is in order [60].)

***Reduction of uniqueness of Problem 5 to the uniqueness of Problem 2*** Let $u$ be the UTM solution of (4.5) [60] and $v$ is a solution of (4.5), satisfying appropriate conditions [60]. Setting $w = (Av + Bv_x + Cv_t) - (Au + Bu_x + Cu_t)$, we have that

$$w_t - w_{xx} = 0, \ w(x,0) = 0 \text{ and } w(0,t) = 0.$$

By Theorem 4, we obtain that $w \equiv 0$, i.e., $A\xi + B\xi_x + C\xi_t \equiv 0$, where we have set $\xi(x,t) = v(x,t) - u(x,t)$. Defining $\zeta(x,t) = \xi(x,t)\exp[Ax/B]$, we obtain the equation $B\zeta_x + C\zeta_t \equiv 0$. Hence $\zeta(x,t) = \phi(\frac{x}{B} - \frac{t}{C})$, for some function $\phi = \phi(\Xi)$, with $\phi(0) = 0$ and $\phi'(0) = 0$. Then $\xi(x,t) = \phi(\frac{x}{B} - \frac{t}{C})\exp[-Ax/B]$ and, therefore,

$$\xi_t - \xi_{xx} = 0 \ \Rightarrow \ \frac{1}{B^2}\phi''(\Xi) + \left(\frac{1}{C} - \frac{2A}{B^2}\right)\phi'(\Xi) - \frac{A^2}{B^2}\phi(\Xi) = 0.$$

Therefore $\phi \equiv 0$, i.e., $v \equiv u$.

Further applications of the new approach to studying uniqueness questions for a wide range of other equations are underway and will appear in subsequent publications.